\documentclass[11pt]{amsart}

\usepackage{amsmath,amsthm, amscd, amssymb, amsfonts}
\usepackage[all]{xy}
\usepackage[dvips, dvipsnames, usenames]{color}

\allowdisplaybreaks
\usepackage{enumitem}

\usepackage{ifthen}

\newcommand{\Dchaintwo}[4]{
	\rule[-3\unitlength]{0pt}{8\unitlength}
	\begin{picture}(14,5)(0,3)
	\put(1,2){\ifthenelse{\equal{#1}{l}}{\circle*{2}}{\circle{2}}}
	\put(2,2){\line(1,0){10}}
	\put(13,2){\ifthenelse{\equal{#1}{r}}{\circle*{2}}{\circle{2}}}
	\put(1,5){\makebox[0pt]{\scriptsize #2}}
	\put(7,4){\makebox[0pt]{\scriptsize #3}}
	\put(13,5){\makebox[0pt]{\scriptsize #4}}
	\end{picture}}

\newcommand{\Dchainthree}[6]{
	\rule[-3\unitlength]{0pt}{8\unitlength}
	\begin{picture}(26,5)(0,3)
	\put(1,2){\ifthenelse{\equal{#1}{l}}{\circle*{2}}{\circle{2}}}
	\put(2,2){\line(1,0){10}}
	\put(13,2){\ifthenelse{\equal{#1}{m}}{\circle*{2}}{\circle{2}}}
	\put(14,2){\line(1,0){10}}
	\put(25,2){\ifthenelse{\equal{#1}{r}}{\circle*{2}}{\circle{2}}}
	\put(1,5){\makebox[0pt]{\scriptsize #2}}
	\put(7,4){\makebox[0pt]{\scriptsize #3}}
	\put(13,5){\makebox[0pt]{\scriptsize #4}}
	\put(19,4){\makebox[0pt]{\scriptsize #5}}
	\put(25,5){\makebox[0pt]{\scriptsize #6}}
	\end{picture}}

\newtheorem{theorem}{Theorem}[section]

\newtheorem{teo intro}{Theorem}

\newtheorem{proposition}[theorem]{Proposition}

\theoremstyle{definition}
\newtheorem{definition}[theorem]{Definition}
\newtheorem{row}{Row}

\theoremstyle{remark}
\newtheorem{remark}[theorem]{Remark}

\def\pf{\begin{proof}}
\def\epf{\end{proof}}

\def\zt{\Z^{\theta}}
\def\z2t{\Z^{2\theta}}

\def\br{\mathfrak{br}}
\def\brj{\mathfrak{brj}}

\def\ufo{\mathfrak{ufo}}

\def\yd{{}^{H}_{H}\mathcal{YD}}

\newcommand\id{\operatorname{id}}

\newcommand\ord{\operatorname{ord}}

\newcommand\co{\operatorname{co}}

\newcommand\ad{\operatorname{ad}}

\def\ot{\otimes}
\def\ra{\rightarrow}

\def\G{\mathbb{G}}
\newcommand\I{\mathbb I}

\def\N{\mathbb{N}}

\def\Z{\mathbb{Z}}

\def\dpn{\widetilde{\mathcal{B}}}

\def\cB{\mathcal{B}}

\def\cH{\mathcal{H}}

\def\cO{\mathcal{O}}
\def\cP{\mathcal{P}}

\def\cU{\mathcal{U}}
\def\cW{\mathcal{W}}
\def\cX{\mathcal{X}}

\newcommand{\g}{\mathfrak g}

\def\n{\mathfrak{n}}
\def\nq{\n_{\bq}}
\def\bq{\mathfrak{q}}

\newcommand{\az}{\mathfrak Z}

\newcommand{\ab}{\mathbf{a}}

\newcommand{\hb}{\mathbf{h}}
\newcommand{\jb}{\mathbf{j}}
\newcommand{\ku}{ \mathbf{k}}

\def\qmb{\mathbf{q}}

\def\xb{x}
\def\yb{\mathbf{y}}

\newcommand{\Ht}{\mathtt H}

\newcommand{\lu}{\mathcal{L}}
\newcommand{\luq}{\lu_{\bq}}

\newcommand{\fO}{\mathfrak O}

\newcommand{\dpndual}{\dpn_{\bq}}

\begin{document}


\title[Lie algebras arising from
Nichols algebras]{A finite-dimensional Lie algebra arising from a Nichols algebra of diagonal type (rank 2)}

\author[Andruskiewitsch; Angiono; Rossi Bertone]
{Nicol\'as Andruskiewitsch, Iv\'an Angiono, Fiorela Rossi Bertone}

\address{FaMAF-CIEM (CONICET), Universidad Nacional de C\'ordoba,
	Medina A\-llen\-de s/n, Ciudad Universitaria, 5000 C\' ordoba, Rep\'
	ublica Argentina.} \email{(andrus|angiono|rossib)@mate.uncor.edu}

\thanks{The work was partially supported by CONICET,
	FONCyT-ANPCyT, Secyt (UNC)}

\begin{abstract}
Let $\cB_\bq$ be a finite-dimensional Nichols algebra of diagonal type corresponding to a matrix $\bq \in \ku^{\theta \times \theta}$.
Let $\lu_{\bq}$ be the Lusztig algebra associated to $\cB_\bq$ \cite{AAR}.
We present $\luq$ as an extension (as braided Hopf algebras) of $\cB_\bq$ by $\az_\bq$ where $\az_\bq$ is isomorphic to the universal enveloping algebra of a Lie algebra $\n_\bq$. We compute the Lie algebra $\n_\bq$ when $\theta = 2$.
\end{abstract}

\maketitle

\section{Introduction}

\subsection{}
Let $\ku$ be a field, algebraically closed
and of characteristic zero.
Let $\theta \in \N$, $\I = \I_{\theta} := \{1, 2, ..., \theta\}$.
Let $\bq=(q_{ij})_{i,j\in\I}$ be a matrix with entries in $\ku^{\times}$, $V$ a vector space with a  basis
$(x_{i})_{i\in\I}$ and $c^{\bq} \in GL(V \ot V)$ be given by
\begin{align*}
c^{\bq}(x_i \ot x_j) &= q_{ij} x_j \ot x_i, & i,j&\in\I.
\end{align*}
Then $(c^{\bq}\otimes \id)(\id\otimes c^{\bq})(c^{\bq}\otimes \id) =
(\id\otimes c^{\bq})(c^{\bq}\otimes \id)(\id\otimes c^{\bq})$, i.e.
$(V, c^{\bq})$ is a braided vector space and the corresponding Nichols algebra $\cB_{\bq} := \cB(V)$ is called of diagonal type.
Recall that $\cB_{\bq}$ is the image of the unique map of braided Hopf algebras $\Omega: T(V) \to T^c(V)$ from the free associative algebra of $V$ to the free associative coalgebra of $V$, such that $\Omega_{\vert V}= \id_V$.
For unexplained terminology and notation, we refer to \cite{AS}.

\medbreak
Remarkably, the explicit classification of all $\bq$ such that $\dim \cB_{\bq} < \infty$ is known  \cite{H-classif RS} (we recall the list when $\theta = 2$ in Table \ref{table:dynkingeneralizados}).
Also, for every $\bq$ in the list of \cite{H-classif RS}, the defining relations
are described in \cite{A-convex, A-presentation}.

\subsection{} Assume that $\dim \cB_{\bq} < \infty$. Two infinite dimensional graded braided Hopf algebras $\widetilde\cB_{\bq}$ and $\lu_{\bq}$ (the Lusztig algebra of $V$) were
introduced and studied in \cite{A-presentation, A-preNichols}, respectively
\cite{AAR}. Indeed,  $\widetilde\cB_{\bq}$ is a pre-Nichols,  and $\lu_{\bq}$  a post-Nichols, algebra of $V$, meaning that $\widetilde\cB_{\bq}$ is intermediate between $T(V)$ and $\cB_{\bq}$, while $\lu_{\bq}$ is intermediate between $\cB_{\bq}$ and $T^c(V)$. This is summarized in the following commutative diagram:
\begin{align*}
\xymatrix{T(V) \ar@/^2pc/^{\Omega}[0,6]
	\ar  @{->}[rrr] \ar  @{->}[1,2] &  & & \cB_{\bq} \ar  @{->}[1,2] \ar  @{->}[rrr] & &  & T^c(V)
	\\
	& &  \widetilde\cB_{\bq} \ar  @{->}_{\pi}[-1,1]   &  & &  \lu_{\bq} \ar  @{->}[-1,1]  &
}
\end{align*}
The algebras $\widetilde\cB_{\bq}$ and $\lu_{\bq}$ are generalizations of the positive parts of the De Concini-Kac-Procesi quantum group, respectively the Lusztig quantum divided powers  algebra.
The distinguished pre-Nichols algebra $\widetilde\cB_{\bq}$ is defined  discarding some of the relations in \cite{A-presentation}, while $\lu_{\bq}$ is  the graded dual of $\widetilde\cB_{\bq}$.

\subsection{}
The following notions are discussed in Section \ref{sec: preliminaries}.
Let $\Delta_+^{\bq}$ be the generalized positive root system of $\cB_\bq$ and let
$\fO_\bq\subset \Delta_+^{\bq}$ be the set of Cartan roots of $\bq$. Let
$x_{\beta}$ be the root vector associated to $\beta\in\Delta_+^{\bq}$, let
$N_\beta = \ord q_{\beta\beta}$ and let $Z_{\bq}$ be the subalgebra of $\dpndual$ generated by $x_{\beta}^{N_{\beta}}$, $\beta \in \fO_{\bq}$.
By \cite[Theorems 4.10, 4.13]{A-preNichols},
$Z_{\bq}$ is a braided normal Hopf subalgebra of $\dpndual$ and $Z_{\bq} = \,^{\co\pi} \dpndual$.  Actually, $Z_{\bq}$ is a true commutative Hopf algebra provided that
\begin{align}\label{condition cart roots}
&q_{\alpha\beta}^{N_\beta}=1, & \forall \alpha,\beta\in\fO_\bq.
\end{align}

Let $\az_\bq$ be the graded dual of $Z_\bq$; under the assumption \eqref{condition cart roots} $\az_\bq$ is a cocommutative Hopf algebra,
hence it is isomorphic to the enveloping algebra $\cU(\nq)$ of the Lie algebra $\nq := \cP(\az_\bq)$.
We show in Section \ref{sec:extensions} that $\luq$ is an extension (as braided Hopf algebras) of $\cB_\bq$ by $\az_\bq$:
\begin{align}\label{eq:extension-braided-lu}
\cB_{\bq} \overset{\pi^*}{\hookrightarrow} \luq \overset{\iota^*}{\twoheadrightarrow} \az_{\bq}.
\end{align}

The main result of this paper is the determination of the Lie algebra $\nq$ when $\theta = 2$ and the generalized Dynkin diagram of $\bq$ is connected.

\begin{theorem}\label{th:main} Assume that $\dim \cB_\bq < \infty$ and $\theta = 2$. Then $\n_\bq$ is either 0 or isomorphic to $\g^+$, where $\g$ is a finite-dimensional semisimple Lie algebra listed in the last column of Table \ref{table:dynkingeneralizados}.
\end{theorem}

Assume that there exists a Cartan matrix $\ab =(a_{ij})$ of finite type,
that becomes symmetric after multiplying with a diagonal $(d_i)$,
and a root of unit $q$ of odd order (and relatively prime to 3 if $\ab$ is of type $G_2$)
such that $q_{ij} = q^{d_ia_{ij}}$ for all $i, j \in \I$. Then \eqref{eq:extension-braided-lu}
encodes the quantum Frobenius homomorphism defined by Lusztig and Theorem \ref{th:main}
is a result from \cite{L-roots of 1}.

\smallbreak
The penultimate column of Table \ref{table:dynkingeneralizados} indicates the type of $\bq$ as established in \cite{AA2}. Thus, we associate  Lie algebras in characteristic zero to some contragredient Lie (super)algebras in positive characteristic.
In a forthcoming paper we shall compute the Lie algebra $\n_\bq$ for $\theta>2$.

\begin{table}
	\setlength{\unitlength}{1mm}
	\begin{tabular}{p{0.6cm}|l|l|l|p{1.35cm}}
		Row & \text{Generalized Dynkin diagrams} & \text{parameters} & \text{Type of} $\cB_\bq$ & $\n_\bq\simeq \g^+$\\
		\hline \hline
		\rule{0mm}{16pt}
		1 & \Dchaintwo{}{$q$}{$q^{-1}$}{$q$} & $q\neq 1$ & \text{Cartan} $A$ & $A_2$ \\
		\hline
		\rule{0mm}{16pt}
		2 & \Dchaintwo{}{$q$}{$q^{-1}$}{$-1$}\ \Dchaintwo{}{$-1$}{$q$}{$-1$} & $q\neq \pm1$ & \text{Super} $A$ & $A_1$\\
		\hline
		\rule{0mm}{16pt}
		3 & \Dchaintwo{}{$q$}{$q^{-2}$}{$q^2$} & $q\neq \pm1$ & \text{Cartan} $B$ & $B_2$\\
		\hline
		\rule{0mm}{16pt}
		4 & \Dchaintwo{}{$q$}{$q^{-2}$}{$-1$} \ \ \Dchaintwo{}{$-q^{-1}$}{$q^2$}{$-1$} & $q\notin \G_4$ & \text{Super} $B$ & $A_1\oplus A_1$\\
		\hline
		\rule{0mm}{16pt}
		5 & \Dchaintwo{}{$\zeta $}{$q^{-1}$}{$q$} \ \ \Dchaintwo{}{$\zeta $}{$\zeta ^{-1}q$}{\ $\zeta q^{-1}$} &
		$\zeta \in \G_3 \not\ni q$ & $\br(2, a)$ & $A_1\oplus A_1$\\
		\hline
		\rule{0mm}{16pt}
		6 & \Dchaintwo{}{$\zeta $}{$-\zeta $}{$-1$} \ \ \Dchaintwo{}{$\zeta ^{-1}$}{$-\zeta ^{-1}$}{$-1$} & $\zeta \in \G_3'$ & \text{Standard} $B$ & $0$\\
		\hline
		\rule{0mm}{16pt}
		7 & \Dchaintwo{}{$-\zeta ^{-2}$}{$-\zeta ^3$}{$-\zeta ^2$}\ \ \Dchaintwo{}{$-\zeta ^{-2}$}{$\zeta ^{-1}$}{$-1$}\ \
		\Dchaintwo{}{$-\zeta ^2$}{$-\zeta $}{$-1$} & $\zeta \in \G_{12}'$ & $\ufo(7)$ & $0$\\
		& \Dchaintwo{}{$-\zeta ^3$}{$\zeta $}{$-1$}\ \ \Dchaintwo{}{$-\zeta ^3$}{$-\zeta ^{-1}$}{$-1$} & & & \\
		\hline
		\rule{0mm}{16pt}
		8 & \Dchaintwo{}{$-\zeta ^2$}{$\zeta $}{$-\zeta ^2$}\ \ \Dchaintwo{}{$-\zeta ^2$}{$\zeta ^3$}{$-1$}\ \ \Dchaintwo{}{$-\zeta ^{-1}$}{$-\zeta ^3$}{$-1$}
		& $\zeta \in \G_{12}'$ & $\ufo(8)$ & $A_1$\\
		\hline
		\rule{0mm}{16pt}
		9 & \Dchaintwo{}{$-\zeta $}{$\zeta ^{-2}$}{$\zeta ^3$}\ \  \Dchaintwo{}{$\zeta ^3$}{$\zeta ^{-1}$}{$-1$}\ \ \Dchaintwo{}{$-\zeta ^2$}{$\zeta $}{$-1$}
		& $\zeta \in \G_9'$ & $\brj(2; 3)$ & $A_1\oplus A_1$\\
		\hline
		\rule{0mm}{16pt}
		10 & \Dchaintwo{}{$q$}{$q^{-3}$}{$q^3$} & $q\notin \G_2 \cup \G_3$ & \text{Cartan} $G_2$ & $G_2$\\
		\hline
		\rule{0mm}{16pt}
		11 & \Dchaintwo{}{$\zeta ^2$}{$\zeta $}{$\zeta ^{-1}$}\ \  \Dchaintwo{}{$\zeta ^2$}{$-\zeta ^{-1}$}{$-1$}\ \  \Dchaintwo{}{$\zeta $}{$-\zeta $}{$-1$}
		& $\zeta \in \G_8'$ & \text{Standard} $G_2$ & $A_1\oplus A_1$\\
		\hline
		\rule{0mm}{16pt}
		12 & \Dchaintwo{}{$\zeta ^6$}{$-\zeta ^{-1}$}{\ \ $-\zeta ^{-4}$}\ \  \Dchaintwo{}{$\zeta ^6$}{$\zeta $}{$\zeta ^{-1}$} & $\zeta \in \G_{24}'$ &
        $\ufo(9)$ & $A_1\oplus A_1$\\
		& \Dchaintwo{}{$-\zeta ^{-4}$}{$\zeta ^5$}{$-1$}\ \ \Dchaintwo{}{$\zeta $}{$\zeta ^{-5}$}{$-1$} & & & \\
		\hline
		\rule{0mm}{16pt}
		13 & \Dchaintwo{}{$\zeta $}{$\zeta ^2$}{$-1$}\ \  \Dchaintwo{}{$-\zeta ^{-2}$}{$\zeta ^{-2}$}{$-1$} & $\zeta \in \G_5'$ & $\brj(2; 5)$ & $B_2$\\
		\hline
		\rule{0mm}{16pt}
		14 & \Dchaintwo{}{$\zeta $}{$\zeta ^{-3}$}{$-1$}\ \ \Dchaintwo{}{$-\zeta $}{$-\zeta ^{-3}$}{$-1$}  & $\zeta \in \G_{20}'$ & $\ufo(10)$ & $A_1\oplus A_1$\\
		& \Dchaintwo{}{$-\zeta ^{-2}$}{$\zeta ^3$}{$-1$}
		\ \ \Dchaintwo{}{$-\zeta ^{-2}$}{$-\zeta ^3$}{$-1$} & & &\\
		\hline
		\rule{0mm}{16pt}
		15 & \Dchaintwo{}{$-\zeta $}{$-\zeta ^{-3}$}{$\zeta ^5$}\ \ \Dchaintwo{}{$\zeta ^3$}{$-\zeta ^4$}{$-\zeta ^{-4}$}
		& $\zeta \in \G_{15}'$ & $\ufo(11)$ & $A_1\oplus A_1$\\
		& \Dchaintwo{}{$\zeta ^5$}{$-\zeta ^{-2}$}{$-1$}\ \ \Dchaintwo{}{$\zeta ^3$}{$-\zeta ^2$}{$-1$} & & &\\
		\hline
		\rule{0mm}{16pt}
		16 & \Dchaintwo{}{$-\zeta $}{$-\zeta ^{-3}$}{$-1$}\ \ \Dchaintwo{}{$-\zeta ^{-2}$}{$-\zeta ^3$}{$-1$} & $\zeta \in \G_7'$ & $\ufo(12)$ & $G_2$\\
		\hline
	\end{tabular}
	\\[2mm] \
	\caption{Lie algebras arising from Dynkin diagrams of rank 2.
	 }\label{table:dynkingeneralizados}
\end{table}

\subsection{}
The paper is organized as follows. We collect the needed preliminary material in Section 2. Section 3 is devoted to  the exactness of \eqref{eq:extension-braided-lu}. The computations of the various $\n_\bq$ is the matter of Section \ref{sec:nq}.
We denote by $\G_N $ the group of $N$-th roots of 1, and by $\G'_N$ its subset of primitive  roots.

\section{Preliminaries}\label{sec: preliminaries}
\subsection{The Nichols algebra, the distinguished-pre-Nichols algebra and the Lusztig algebra}
Let $\bq$ be as in the Introduction and let $(V,c^{\bq})$ be the corresponding braided vector space of diagonal type. We  assume from now on that $\cB_\bq$ is finite-dimensional. Let $(\alpha_j)_{j\in \I}$ be the canonical basis of $\zt$.
Let $\qmb:\zt\times\zt\to\ku^\times$ be the $\Z$-bilinear form associated to the matrix $\bq$, i.e. $\qmb(\alpha_j,\alpha_k)=q_{jk}$ for all $j,k\in\I$.
If $\alpha,\beta  \in \zt$, we set
$ q_{\alpha\beta} = \qmb(\alpha,\beta)$.
Consider the matrix $(c_{ij}^{\bq})_{i,j\in \I}$, $c_{ij}\in\Z$
defined by $c_{ii}^{\bq} = 2$,
\begin{align}\label{eq:defcij}
c_{ij}^{\bq}&:= -\min \left\{ n \in \N_0: (n+1)_{q_{ii}}
(1-q_{ii}^n q_{ij}q_{ji} )=0 \right\},  & i & \neq j.
\end{align}
This is well-defined by \cite{R}.
Let $i\in \I$. We recall the following definitions:

\begin{itemize}[leftmargin=*]\renewcommand{\labelitemi}{$\diamond$}
	\item The reflection
	$s_i^{\bq}\in GL(\Z^\theta)$, given by $s_i^{\bq}(\alpha_j)=\alpha_j-c_{ij}^{\bq}\alpha_i$, $j\in \I$.
	
	\medbreak
	\item  The matrix  $\rho_i(\bq)$, given by
	$ \rho_i(\bq)_{jk}= \qmb(s_i^{\bq}(\alpha_j),s_i^{\bq}(\alpha_k))$, $j, k \in \I$.
	
		\medbreak
		\item The braided vector space $\rho_i(V)$ of diagonal type with matrix  $\rho_i(\bq)$.
\end{itemize}
A basic result is that $\cB_{\bq} \simeq \cB_{\rho_i(\bq)}$, at least as graded vector spaces.

The algebras $T(V)$ and $\cB_{\bq}$ are $\zt$-graded by $\deg x_i = \alpha_i$, $i\in \I$.
Let $\Delta_+^{\bq}$ be the set of $\zt$-degrees of the generators of a PBW-basis of $\cB_\bq$, counted with multiplicities \cite{H-Weyl grp}.
The elements of $\Delta_+^{\bq}$ are called (positive) roots.
Let $\Delta^{\bq} = \Delta_+^{\bq} \cup -\Delta_+^{\bq}$. Let
\begin{align*}
\cX := \{\rho_{j_1} \dots \rho_{j_N}(\bq): j_1, \dots, j_N \in \I, N \in \N \}.
\end{align*}
Then the  generalized  root system of $\bq$ is the fibration $\Delta \to  \cX$, where the fiber of $\rho_{j_1} \dots \rho_{j_N}(\bq)$ is $\Delta^{\rho_{j_1} \dots \rho_{j_N}(\bq)}$. The Weyl groupoid of $\cB_\bq$ is a groupoid, denoted
$\cW_{\bq}$, that acts on this fibration, generalizing the classical Weyl group, see \cite{H-Weyl grp}. We know from \emph{loc. cit.} that $\cW_{\bq}$ is finite (and this characterizes finite-dimensional Nichols algebras of diagonal type).

Here is a useful description of $\Delta_+^{\bq}$.
Let $w \in \cW_{\bq}$ be an element of maximal length. We fix a reduced expression $w=\sigma_{i_1}^{\bq} \sigma_{i_2}\cdots \sigma_{i_M}$.
For $1\leq k\leq M$ set
\begin{align} \label{eq:betak}
\beta_k &= s_{i_1}^{\bq}\cdots s_{i_{k-1}}(\alpha_{i_k}),
\end{align}
Then $\Delta_+^{\bq}=\{\beta_k|1\leq k\leq M\}$  \cite[Prop. 2.12]{CH1};
in particular $\vert \Delta_+^{\bq} \vert = M$.

\smallbreak The notion of Cartan root is instrumental for the definitions of $\dpn_{\bq}$ and $\luq$.
First, following \cite{A-preNichols} we say that $i\in\I$ is a \emph{Cartan vertex}  of $\bq$  if
\begin{align}\label{eq:cartan-vertex}
q_{ij}q_{ji} &= q_{ii}^{c_{ij}^{\bq}}, & \text{for all } j \neq i,
\end{align}
Then the set of \emph{Cartan roots} of $\bq$ is
\begin{align*}
\fO_{\bq} &= \{s_{i_1}^{\bq} s_{i_2} \dots s_{i_k}(\alpha_i) \in \Delta_+^{\bq}:
i\in \I  \text{ is a Cartan vertex of } \rho_{i_k} \dots \rho_{i_2}\rho_{i_1}(\bq) \}.
\end{align*}

Given a positive root $\beta\in\Delta_+^{\bq}$, there is an associated
root vector $x_{\beta} \in \cB_{\bq}$ defined via the so-called Lusztig isomorphisms \cite{H-isom}.
Set $N_\beta = \ord q_{\beta\beta} \in\N$,  $\beta \in \Delta_+^{\bq}$.
Also, for $\hb=(h_1,\dots,h_M)\in\N_0^{M}$ we write
\begin{align*}
\xb^{\hb}=x_{\beta_M}^{h_M}x_{\beta_{M-1}}^{h_{M-1}} \cdots
x_{\beta_1}^{h_1}.
\end{align*}
Let $\widetilde N_k = \begin{cases} N_{\beta_k} &\mbox{ if }\beta_k\notin\cO_{\bq},
\\ \infty  &\mbox{ if }\beta_k\in\cO_{\bq}. \end{cases}$.
For simplicity, we introduce
\begin{align}\label{eq:ht}
\Ht = \{\hb\in\N_0^M: \, 0\leq h_k < \widetilde N_k, \text{ for all } k\in \I_M \}.
\end{align}
By \cite[Theorem 3.6]{A-preNichols}
the set $\{ \xb^{\hb}\, | \,  \hb\in\Ht \}$ is a basis of $\dpn_{\bq}$.

\smallbreak
As said in the Introduction, the Lusztig algebra associated to $\cB_\bq$ is the braided Hopf algebra $\luq$ which is the graded dual of $\dpndual$.
Thus, it comes equipped with
a bilinear form $\langle \, , \, \rangle: \dpndual \times \luq \rightarrow \ku$, 
which satisfies for all $x,x' \in \dpndual$, $y,y' \in \luq$
\begin{align*}
&\langle y, xx'\rangle = \langle y^{(2)},x\rangle\langle y^{(1)},x' \rangle \qquad \mbox{ and }
\qquad \langle yy', x\rangle = \langle y,x^{(2)}\rangle\langle y',x^{(1)} \rangle.
\end{align*}
If $\hb\in\Ht$, then define $\yb_{\hb} \in \luq$  by $\langle \yb_{\hb},\xb^{\jb}\rangle = \delta_{\hb, \jb}$, $\jb \in \Ht$.
Let $(\hb_k)_{k\in \I_M}$ denote the canonical basis of $\Z^M$.
If $k\in \I_M$ and $\beta=\beta_k \in \Delta_+^{\bq}$, then we denote the element $\yb_{n\hb_k}$ by $y_{\beta}^{(n)}$.
Then the algebra $\luq$ is generated by
\begin{align*}
\{ y_{\alpha}: \alpha\in \Pi_{\bq}\}  \cup  \{y_{\alpha}^{(N_\alpha)}: \alpha\in \fO_{\bq}, \,  x_{\alpha}^{N_{\alpha}}\in\mathcal{P}(\dpndual) \},
\end{align*}
by \cite{AAR}. Moreover, by \cite[4.6]{AAR}, the  following set is a basis of $\luq$:
$$\{ y_{\beta_1}^{(h_1)} \cdots y_{\beta_M}^{(h_M)} |\, (h_1,\dots,h_M)\in\Ht \}.$$

\subsection{Lyndon words, convex order and PBW-basis}
For the computations in Section \ref{sec:nq} we need some preliminaries on Kharchenko's PBW-basis.
Let $(V,\bq)$ be as above and let $\mathbb{X}$ be the set of words with letters in $X=\{x_1,\dots,x_\theta\}$ (our fixed basis of $V$); the empty word is $1$ and for $u\in\mathbb{X}$ we write $\ell(u)$ the length of $u$.  We can identify $\ku \mathbb{X}$ with $T(V)$.

\begin{definition} Consider the lexicographic order in  $\mathbb{X}$.
We say that $u \in \mathbb{X} -\{ 1\}$ is a \emph{Lyndon word} if for every decomposition $u=vw$, $v,w \in\mathbb{X} - \left\{ 1 \right\}$, then $u<w$.
We denote by $L$ the set of all Lyndon words.
\end{definition}
	
A well-known theorem, due to Lyndon, established that any word $u \in \mathbb{X}$ admits a unique decomposition, named \emph{Lyndon decomposition}, as a non-increasing product of Lyndon words:
\begin{equation}\label{eq:descly}
u=l_1l_2\dots  l_r, \qquad l_i \in L, l_r \leq \dots \leq l_1.
\end{equation}
Also, each $l_i\in L$ in \eqref{eq:descly} is called a
\emph{Lyndon letter} of $u$.

Now each  $u \in L - X$ admits at least one decomposition $u=v_1v_2$ with $v_1,v_2\in L$.
Then the \emph{Shirshov decomposition} of $u$ is the decomposition $u=u_1u_2$, $u_1,u_2\in L$, such
that $u_2$ is the smallest end of $u$ between all possible decompositions of this form.

\medskip
For any braided vector space $V$, the \emph{braided bracket} of $x,y\in T(V)$ is
\begin{equation}\label{eq:braidedcommutator}
[x,y]_c := \text{multiplication } \circ \left( \id - c \right) \left( x \ot y \right).
\end{equation}

Using the identification $T(V)=\ku \mathbb{X}$ and the decompositions described above, we can define a $\ku$-linear
endomorphism $\left[ - \right]_c$ of $T(V)$ as follows:
$$ \left[ u \right]_c := \begin{cases} u,& \text{if } u = 1 \text{ or }u \in X;\\
[\left[ v \right]_c, \left[ w \right]_c]_c,  & \text{if } u \in L-X, \ u=vw \text{\small \ its Shirshov decomposition};\\
\left[u_1\right]_c \dots \left[u_t\right]_c,& \text{if } u\in \mathbb{X}-L, u=u_1\dots u_t \text{\small \  its Lyndon decomposition}.
\end{cases}
$$
We will describe PBW-bases using this endomorphism.

\begin{definition} For $l \in L$, the element $\left[l\right]_c$ is the corresponding \emph{hyperletter}. A word written in hyperletters is an \emph{hyperword}; a \emph{monotone hyperword} is an hyperword $W=\left[u_1\right]_c^{k_1}\dots\left[u_m\right]_c^{k_m}$
such that $u_1>\dots >u_m$.
\end{definition}

Consider now a different order on $\mathbb{X}$, called \emph{deg-lex order} \cite{Kh}: For each pair $u,v \in \mathbb{X}$, we have that $u \succ v$ if
$\ell(u)<\ell(v)$, or $\ell(u)=\ell(v)$ and $u>v$ for the
lexicographical order. This order is total, the empty word $1$ is the maximal
element and it is invariant by left and right multiplication.

\smallbreak
Let $I$ be a Hopf ideal of $T(V)$ and $R=T(V)/I$. Let $\pi: T(V) \rightarrow R$ be the canonical
projection. We set:
$$G_I:= \left\{ u \in \mathbb{X}: u \notin \\ \ku \mathbb{X}_{\succ u}+I  \right\}.$$
Thus, if $u \in G_I$ and $u=vw$, then $v,w \in G_I$. So, each $u \in G_I$ is a non-increasing product of Lyndon words of $G_I$.

Let $S_I:=G_I\cap L$ and let $h_I:S_I\to\left\{2,3,\dots \right\}\cup \{\infty\}$ be defined by:
\begin{equation}\label{defaltura}
    h_I(u):= \min \left\{ t \in \N : u^t  \in \ku \mathbb{X}_{\succ u^t} + I \right\}.
\end{equation}

\begin{theorem}\label{thm:base PBW Kharchenko} \cite{Kh} The following set is a PBW-basis of $R = T(V)/I$:
$$ \{ \left[u_1\right]_c^{k_1}\dots\left[u_m\right]_c^{k_m}: \, m\in\N_0, u_1> \ldots >u_m, u_i\in S_I, 0<k_i<h_I(u_i) \}. \qed$$
\end{theorem}

We refer to this base as \emph{Kharchenko's PBW-basis} of $T(V)/I$ (it depends on the order of $X$).

\medskip

\begin{definition}\cite[2.6]{A-convex}
Let $\Delta_\bq^+$ be as above and let $<$ be a total order on $\Delta_\bq^+$. We say that the order is \emph{convex} if for each $\alpha,\beta\in\Delta_\bq^+$ such that $\alpha<\beta$ and $\alpha+\beta\in\Delta_\bq^+$, then $\alpha<\alpha+\beta<\beta$.
The order is called	\emph{strongly convex} if for each ordered subset $\alpha_1\leq\alpha_2\leq\dots\leq\alpha_k$ of elements of $\Delta_\bq^+$ such that $\alpha=\sum_i \alpha_i\in\Delta_\bq^+$, then $\alpha_1<\alpha<\alpha_k$.
\end{definition}

\begin{theorem}\cite[2.11]{A-convex}
The following statements are equivalent:
\begin{itemize}[leftmargin=*]
	\item The order is convex.
	\item The order is strongly convex.
	\item The order arises from a reduced expression of a longest element $w \in \cW_{\bq}$, cf. \eqref{eq:betak}. \qed
	\end{itemize}
\end{theorem}

Now, we have two PBW-basis of $\cB_q$ (and correspondingly of $\dpn_{\bq}$), namely Kharchenko's PBW-basis and the PBW-basis defined from a reduced expression of a longest element of the Weyl groupoid. But both basis are reconciled by
\cite[Theorem 4.12]{AY}, thanks to \cite[2.14]{A-convex}.
Indeed, each generator  of Kharchenko's PBW-basis  is a multiple scalar of a generator of the secondly mentioned PBW-basis. So, for ease of calculations, in the rest of this work we shall use the Kharchenko generators.

\smallbreak
The following proposition is used to compute the hyperword $[l_\beta]_c$ associated to a root $\beta\in\Delta_\bq^+$:
\begin{proposition}\cite[2.17]{A-convex}
	For $\beta\in\Delta_\bq^+$,
	$$
	l_\beta=\begin{cases}
	x_{\alpha_i}, & \mbox{ if } \beta=\alpha_i, \, i\in\I; \\
	\max\{l_{\delta_1}l_{\delta_2}: \, \delta_1,\delta_2\in\Delta_\bq^+, \delta_1+\delta_2=\beta, l_{\delta_1}<l_{\delta_2} \} , & \mbox{ if } \beta\neq\alpha_i, \, i\in\I.\qed
	\end{cases}
	$$
\end{proposition}

We give a list of the hyperwords appearing in the next section:
$$
\begin{array}{ccc}
\text{Root} & \text{Hyperword} & \text{Notation} \\ \hline
\alpha_i & x_i & x_i \\
n\alpha_1+\alpha_2 & (\ad_cx_1)^nx_2 & x_{1\dots 12} \\
\alpha_1+2\alpha_2 & [x_{\alpha_1+\alpha_2},x_2]_c & [x_{12},x_2]_c \\
3\alpha_1+2\alpha_2 & [x_{2\alpha_1+\alpha_2},x_{\alpha_1+\alpha_2}]_c & [x_{112},x_{12}]_c \\
4\alpha_1+3\alpha_2 & [x_{3\alpha_1+2\alpha_2},x_{\alpha_1+\alpha_2}]_c & [[x_{112},x_{12}]_c,x_{12}]_c \\
5\alpha_1+3\alpha_2 & [x_{2\alpha_1+\alpha_2},x_{3\alpha_1+2\alpha_2}]_c & [x_{112},[x_{112},x_{12}]_c]_c
\end{array}
$$
We use an analogous notation for the elements of $\luq$: for example we write $y_{112,12}$ when we refer to the element of $\luq$ which corresponds to $[x_{112},x_{12}]_c$.

\section{Extensions of braided Hopf algebras}\label{sec:extensions}

We recall the definition of braided Hopf algebra extensions given in \cite{AN}; we refer to \cite{BD, GG} for  more general definitions. Below we denote by $\underline{\Delta}$ the coproduct of a braided Hopf algebra $A$ and by $A^+$ the kernel of the counit.

\smallbreak
 First, if $\pi: C \to B$
is a morphism of Hopf algebras in $\yd$, then we set
\begin{align*}
C^{\,\co\pi} &=\{c\in C \,| \, (\id\ot\pi)\underline{\Delta}(c)=c\ot 1\}, \\
\,^{\co\pi}C &=\{c\in C \,| \, (\pi\ot\id)\underline{\Delta}(c)=1\ot c\}.
\end{align*}
\begin{definition}\cite[\S 2.5]{AN}\label{def: ext bha}
Let $H$ be a Hopf algebra. A sequence of morphisms of Hopf algebras in $\yd$
	\begin{align}\label{eq:def-exseq}
	\ku \rightarrow A \overset{\iota}{\to} C \overset{\pi}{\to} B \rightarrow \ku
	\end{align}
is an \emph{extension of braided Hopf algebras} if
	\begin{itemize}
		\item[(i)] $\iota$ is injective,
		\item[(ii)] $\pi$ is surjective,
		\item[(iii)] $\ker \pi = C\iota(A^+)$ and
		\item[(iv)] $A=C^{\,\co\pi}$, or equivalently $A=\,^{\co\pi}C$.
	\end{itemize}
\end{definition}

For simplicity, we shall write $A \overset{\iota}{\hookrightarrow} C \overset{\pi}{\twoheadrightarrow} B$ instead of \eqref{eq:def-exseq}.

\smallbreak
This Definition applies in our context: recall that  $\cB_{\bq} \simeq\dpndual /\langle x_{\beta}^{N_{\beta}}, \,\beta \in \fO_{\bq}\rangle$.
Let $Z_{\bq}$ be the subalgebra of $\dpndual$ generated by $x_{\beta}^{N_{\beta}}$, $\beta \in \fO_{\bq}$. Then

\begin{itemize}[leftmargin=*]\renewcommand{\labelitemi}{$\circ$}
	\item The inclusion $\iota: Z_{\bq} \to \dpndual$ is injective and the projection $\pi: \dpndual \to \cB_{\bq}$ is surjective.
	
	\smallbreak
	\item  \cite[Theorem 4.10]{A-preNichols} $Z_{\bq}$ is a \emph{normal} Hopf subalgebra of $\dpndual$; since $\ker \pi$ is the two-sided ideal generated by $\iota(Z_{\bq}^+)$,  $\ker \pi = \dpndual\iota(Z_{\bq}^+)$.
	
	\smallbreak
	\item \cite[Theorem 4.13]{A-preNichols} $Z_{\bq} = \,^{\co\pi} \dpndual$.
\end{itemize}

Hence we have an extension of braided Hopf algebras
\begin{align}\label{eq:extension-braided-dpn}
 Z_{\bq} \overset{\iota}{\hookrightarrow} \dpn_{\bq} \overset{\pi}{\twoheadrightarrow} \cB_{\bq}.
\end{align}

The morphisms $\iota$ and $\pi$ are graded. Thus, taking graded duals, we obtain a new sequence of morphisms of braided Hopf algebras
\begin{align}\tag{\ref{eq:extension-braided-lu}}
\cB_{\bq} \overset{\pi^*}{\hookrightarrow} \luq \overset{\iota^*}{\twoheadrightarrow} \az_{\bq}.
\end{align}
\begin{proposition}
	The sequence \eqref{eq:extension-braided-lu} is an extension of braided Hopf algebras.
\end{proposition}

\pf The argument of \cite[3.3.1]{A-canad} can be adapted to the present situation, or more generally to extensions of braided Hopf algebras that are graded with finite-dimensional homogeneous components.
The map $\pi^*:\cB_{\bq}\rightarrow\luq$ is injective because $\cB_{\bq} \simeq \cB_{\bq}^*$;
$\iota^*: \luq \overset{\iota^*}\to \az_{\bq}$ is surjective being the transpose of a graded monomorphism between two locally finite graded vector spaces.
Now, since $Z_{\bq} = \,^{\co\pi} \dpndual=   \dpndual^{\,\co\pi}$, we have
\begin{align}\label{eq:coinv-exam}
\ker\iota^*=\luq\cB_\bq^+=\cB_\bq^+\luq.
\end{align}
Similarly $\luq^{\,\co\iota^*} = \cB_{\bq}^*$ because  $\ker \pi^{\bot} = \cB_{\bq}$.
\epf

From now on, we assume  the  condition \eqref{condition cart roots} on the matrix $\bq$ mentioned in the Introduction, that is
\begin{align*}
&q_{\alpha\beta}^{N_\beta}=1, & \forall \alpha,\beta\in\fO_\bq.
\end{align*}
The following result is our basic tool to compute the Lie algebra $\n_\bq$.

\begin{theorem}\label{th:nq}	
The braided Hopf algebra $\az_\bq$ is an usual Hopf algebra, isomorphic to
the universal enveloping algebra of the Lie algebra $\n_\bq=\cP(\az_{\bq})$.
The elements $\xi_\beta:=\iota^*(y_\beta^{(N_\beta)})$, $\beta\in\fO_\bq$, form a basis of $\n_\bq$.
\end{theorem}
\pf
Let $A_\bq$ be the subspace of $\luq$ generated by the ordered monomials $y_{\beta_{i_1}}^{(r_1N_{\beta_{i_1}})} \dots y_{\beta_{i_k}}^{(r_kN_{\beta_{i_k}})}$ where $\beta_{i_1}<\dots<\beta_{i_k}$ are all the Cartan roots of $\cB_\bq$ and $r_1,\dots,r_k\in\N_0$.
We claim that  the restriction of the multiplication $\mu:\cB_\bq\ot A_\bq\ra\luq$ is an isomorphism of vector spaces.
Indeed, $\mu$ is surjective by the commuting relations in $\luq$.
Also, the Hilbert series of $\luq$, $\cB_\bq$ and $A_\bq$ are respectively:
\begin{align*}
\cH_{\luq} &=\prod_{\beta_k\in\fO_{\bq}} \frac{1}{1 - T^{\deg \beta}}.\prod_{\beta_k\notin\fO_{\bq}} \frac{1 - T^{N_\beta\deg \beta}}{1 - T^{\deg \beta}};\\
\cH_{\cB_\bq} &=\prod_{\beta_k\in\Delta^+_{\bq}} \frac{1 - T^{N_\beta\deg \beta}}{1 - T^{\deg \beta}};\\
\cH_{A_\bq} &=\prod_{\beta_k\in\fO_{\bq}} \frac{1}{1 - T^{N_\beta\deg \beta}}.
\end{align*}
Since the multiplication is graded and $\cH_{\luq}=\cH_{\cB_\bq}\cH_{A_\bq}$, $\mu$ is injective.
The claim follows and we have
\begin{align}\label{eq:decomp}
\luq = A_\bq \oplus \cB_\bq^+A_\bq.
\end{align}

We next claim that $\iota^*: A_\bq \to \az_\bq$ is an isomorphism of vector spaces.
Indeed, by \eqref{eq:coinv-exam}, $\ker\iota^*=\cB_\bq^+\luq= \cB_\bq^+(\cB_\bq A_\bq)=\cB_\bq^+ A_\bq$. By \eqref{eq:decomp}, the claim follows.

By \eqref{condition cart roots}, $Z_\bq$ is a commutative Hopf algebra, see \cite{A-preNichols}; hence $\az_\bq$ is a cocommutative Hopf algebra. Now the elements  $\xi_\beta:=\iota^*(y_\beta^{(N_\beta)})$, , $\beta\in\fO_\bq$, are primitive, i.e. belong to $\n_\bq=\cP(\az_{\bq})$. The monomials
$\xi_{\beta_{i_1}}^{r_1}\dots \xi_{\beta_{i_k}}^{r_k}$, $\beta_{i_1}<\dots<\beta_{i_k}\in\fO_\bq$, $r_1,\dots,r_k\in\N_0$
form a basis of $\az_\bq$, hence
\begin{align*}
\az_\bq = \ku\langle \xi_{\beta}: \beta\in\fO_\bq\rangle \subseteq \cU(\n_\bq)  \subseteq \az_\bq.
\end{align*}
We conclude that $(\xi_{\beta})_{\beta\in\fO_\bq}$ is a basis of $\n_\bq$ and that
$\az_\bq = \cU(\n_\bq)$.
\epf

\section{Proof of Theorem \ref{th:main}}\label{sec:nq}

In this section we consider all indecomposable matrices $\bq$ of rank 2 whose associated Nichols algebra $\cB_\bq$
is finite-dimensional; these are classified in \cite{H-classif RS} and we recall their diagrams in Table \ref{table:dynkingeneralizados}.
For each $\bq$ we obtain an isomorphism between $\az_\bq$ and $\cU(\g^+)$, the universal enveloping algebra of the positive part of $\g$.
Here $\g$ is the semisimple Lie algebra of the last column of Table \ref{table:dynkingeneralizados},  
with Cartan matrix $A=(a_{ij})_{1\le i,j\le 2}$. By simplicity we denote $\g$ by its type, e.g. $\g = A_2$.

We recall that we assume \eqref{condition cart roots} and that $\xi_\beta= \iota^*(y_\beta^{(N_\beta)}) \in \az_\bq$. 
Thus,
\begin{align*}
[\xi_\alpha, \xi_\beta]_c &= \xi_\alpha\xi_\beta-\xi_\beta\xi_\alpha =[\xi_\alpha, \xi_\beta],
& \mbox{for all } & \alpha,\beta\in\fO_\bq.
\end{align*}

The strategy to prove the isomorphism $\mathfrak{F}:\cU(\g^+) \to \az_\bq$ is the following:
\begin{enumerate}[leftmargin=*]
  \item If $\fO_\bq=\emptyset$, then $\g^+=0$. If $|\fO_\bq|=1$, then $\g=\mathfrak{sl}_2$, i.e. of type $A_1$.
  \item If $|\fO_\bq|=2$, then $\g$ is of type $A_1\oplus A_1$. Indeed, let $\fO_\bq=\{\alpha,\beta\}$. As $\az_\bq$ is
  $\N_0^{\theta}$-graded, $[\xi_\alpha, \xi_\beta]\in\n_\bq$ has degree $N_\alpha \alpha+N_\beta \beta$. Thus
  $[\xi_\alpha, \xi_\beta]=0$.
  \item Now assume that $|\fO_\bq|>2$. We recall that $\az_\bq$ is generated by
  $$ \{ \xi_\beta | x_\beta^{N_\beta} \mbox{ is a  primitive element of }\dpndual \} .$$
  We compute the coproduct of all $x_\beta^{N_\beta}$ in $\dpndual$, $\beta\in\fO_\bq$, 
  using that $\underline{\Delta}$ is a graded map
and $Z_\bq$ is a Hopf subalgebra of $\dpndual$. In all cases we get two primitive elements $x_{\beta_1}^{N_{\beta_1}}$ and $x_{\beta_2}^{N_{\beta_2}}$, thus $\az_\bq$ is generated by $\xi_{\beta_1}$ and $\xi_{\beta_2}$.
  \item Using the coproduct again, we check that 
 \begin{align}\label{eq:serre}
 (\ad \xi_{\beta_i})^{1-a_{ij}} \xi_{\beta_j} &=0, & 1&\le i\neq j\le 2.
 \end{align}  
 To prove \eqref{eq:serre}, it is enough to observe that $\n_\bq$ has a trivial component of degree $N_{\beta_i}(1-a_{ij})\beta_i+N_{\beta_j} \beta_j$.
  Now \eqref{eq:serre} implies that there exists a surjective map of Hopf algebras $\mathfrak{F}:\cU(\g^+) \twoheadrightarrow \az_\bq$ such that $e_i\mapsto\xi_{\beta_i}$.
  \item To prove that $\mathfrak{F}$ is an isomorphism, it suffices to see that the restriction $\g^+ \overset{\ast}{\to} \n_\bq$ is an isomorphism; but in each case we see that $\ast$ is surjective, and  $\dim \g^+=\dim \n_\bq = \vert \fO_\bq \vert$.
\end{enumerate}

\smallskip
We refer to \cite{A-standard, AAY, A-unident} for the presentation, root system and Cartan roots of braidings of standard, super and unidentified type respectively.

\begin{row}
Let $q\in\G_N'$, $N \geq 2$. The diagram
\setlength{\unitlength}{1mm} \Dchaintwo{}{$q$}{$q^{-1}$}{$q$}
corresponds to a braiding of Cartan type $A_2$ whose set of positive roots is $\Delta^+_\bq=\{ \alpha_1, \alpha_1+\alpha_2, \alpha_2\}$.
In this case $\fO_\bq=\Delta^+_\bq$ and $N_\beta=N$ for all $\beta\in\fO_\bq$. By hypothesis, $q_{12}^N=q_{21}^N=1$.
The elements $x_1,x_2\in\dpn_\bq$ are primitive and
$$\underline{\Delta}(x_{12})=x_{12}\ot 1+ 1\ot x_{12}+(1-q^{-1}) x_1\ot x_2. $$
Then the coproducts of the elements $x_1^N,x_{12}^N,x_2^N\in\dpn_\bq$ are:
	\begin{align*}
	&\underline{\Delta}(x_1^N)=x_1^N\ot 1+ 1\ot x_1^N; \qquad
	\underline{\Delta}(x_2^N)=x_2^N\ot 1+ 1\ot x_2^N; \\
	&\underline{\Delta}(x_{12}^N)=x_{12}^N\ot 1+ 1\ot x_{12}^N
	+(1-q^{-1})^Nq_{21}^{\frac{N(N-1)}{2}} x_1^N\ot x_2^N.
	\end{align*}
As $[\xi_2,\xi_{12}]$, $[\xi_1,\xi_{12}]\in\n_\bq$ have degree $N\alpha_1+2N\alpha_2$, respectively $2N\alpha_1+N\alpha_2$, and
the components of these degrees of $\n_\bq$ are trivial, we have
\begin{align*}
	&[\xi_2,\xi_{12}]=[\xi_1,\xi_{12}]=0.
\end{align*}
Again by degree considerations, there exists $c\in\ku$ such that $[\xi_2,\xi_1]=c\xi_{12}$.
By the duality between $\az_\bq$ and $Z_\bq$ we have that
\begin{align*}
	&[\xi_2,\xi_1]= (1-q^{-1})^Nq_{21}^{\frac{N(N-1)}{2}} \xi_{12}.	
\end{align*}
Then there exists a morphism of algebras $\mathfrak{F}: \cU(A_2^+)\ra \az_\bq$ given by
	$$ e_1\mapsto \xi_1, \qquad e_2\mapsto \xi_2.$$
This morphism takes a basis of $A_2^+$ to a basis of $\n_\bq$, so $\az_\bq\simeq \cU(A_2^+)$.
\end{row}

\begin{row}
Let $q\in\G_N'$, $N \geq 3$.
These diagrams correspond to braidings of super type $A$ with positive roots
$\Delta^+_\bq=\{ \alpha_1, \alpha_1+\alpha_2, \alpha_2\}$.
	
The first diagram is \setlength{\unitlength}{1mm} \Dchaintwo{}{$q$}{$q^{-1}$}{$-1$}.	 
In this case the unique Cartan root is $\alpha_1$ with $N_{\alpha_1}=N$.
The element $x_1^N\in\dpn_\bq$ is primitive and $\az_\bq$ is generated by $\xi_1$.
Hence $\az_\bq\simeq \cU(A_1^+)$.
	
The second diagram gives a similar situation, since $\fO_\bq=\{\alpha_1+\alpha_2\}$.
\end{row}

\begin{row}
Let $q\in\G_N'$, $N \geq 3$. The diagram \setlength{\unitlength}{1mm}\Dchaintwo{}{$q$}{$q^{-2}$}{$q^2$}
corresponds to a braiding of Cartan type $B_2$ with $\Delta^+_\bq=\{ \alpha_1, 2\alpha_1+\alpha_2, \alpha_1+\alpha_2, \alpha_2\}$.
In this case $\fO_\bq=\Delta_\bq^+$. The coproducts of the generators of $\dpn_\bq$ are:
\begin{align*}
\underline{\Delta}(x_1)=&x_1\ot 1+ 1\ot x_1; \qquad
\underline{\Delta}(x_2)=x_2\ot 1+ 1\ot x_2; \\
\underline{\Delta}(x_{12})=&x_{12}\ot 1+ 1\ot x_{12}
+(1-q^{-2})\, x_1\ot x_2; \\
\underline{\Delta}(x_{112})=&x_{112}\ot 1+ 1\ot x_{112}
+(1-q^{-1})(1-q^{-2})\, x_1^{2}\ot x_2 \\&+q(1-q^{-2})\, x_1\ot x_{12}.
\end{align*}
	
We have two different cases depending on the parity of $N$.

\begin{enumerate}
  \item If $N$ is odd, then $N_\beta=N$ for all $\beta\in\Delta_\bq^+$. In this case,
	\begin{align*}
	\underline{\Delta}(x_1^N)=&x_1^N\ot 1+ 1\ot x_1^N; \qquad
	\underline{\Delta}(x_2^N)=x_2^N\ot 1+ 1\ot x_2^N; \\
	\underline{\Delta}(x_{12}^N)=&x_{12}^N\ot 1+ 1\ot x_{12}^N
	+(1-q^{-2})^N x_1^N\ot x_2^N; \\
	\underline{\Delta}(x_{112}^N)=&x_{112}^N\ot 1+ 1\ot x_{112}^N
	+(1-q^{-1})^N(1-q^{-2})^N x_1^{2N}\ot x_2^N \\&
	+C\, x_1^N\ot x_{12}^N,
	\end{align*}
\end{enumerate}
for some $C\in\ku$. Hence, in $\az_\bq$ we have the relations
	\begin{align*}
	&[\xi_1,\xi_2]=(1-q^{-2})^N \xi_{12}; \\
	&[\xi_{12},\xi_{1}]=C\, \xi_{112};\\
	&[\xi_1, \xi_2]_c= (1-q^{-1})^N(1-q^{-2})^N \xi_{112} +(1-q^{-2})^N\xi_1\xi_{12};\\
	&[\xi_1,\xi_{112}] =[\xi_2,\xi_{12}] =0.
	\end{align*}
Thus there exists an algebra map $\mathfrak{F}: \cU(B_2^+)\ra \az_\bq$ given by
$ e_1\mapsto \xi_1$, $e_2\mapsto \xi_2$.
Moreover, $\mathfrak{F}$ is an isomorphism, and so $\az_\bq\simeq \cU(B_2^+)$.
Using the relations of $\cU(B_2^+)$ we check that $C=2(1-q^{-1})^N(1-q^{-2})^N$.
\medbreak

\begin{enumerate}\item[(2)] If $N=2M>2$, then $N_{\alpha_1}=N_{\alpha_1+\alpha_2}=N$ and $N_{2\alpha_1+\alpha_2}=N_{\alpha_2}=M$.
In this case we have
	\begin{align*}
	\underline{\Delta}(x_1^N)=&x_1^N\ot 1+ 1\ot x_1^N; \qquad
	\underline{\Delta}(x_2^M)=x_2^M\ot 1+ 1\ot x_2^M; \\
	\underline{\Delta}(x_{12}^N)=&x_{12}^N\ot 1+ 1\ot x_{12}^N
	+(1-q^{-2})^N q_{21}^{M(N-1)}\, x_1^N\ot x_2^{2M} \\&+ (1-q^{-2})^M q_{21}^{M^2} x_{112}^M\ot x_2^M; \\
	\underline{\Delta}(x_{112}^M)=&x_{112}^M\ot 1+ 1\ot x_{112}^M
	+(1-q^{-1})^M(1-q^{-2})^M q_{21}^{M(M-1)} x_1^{N}\ot x_2^M.
	\end{align*}
\end{enumerate}
Hence, the following relations hold in $\az_\bq$:
	\begin{align*}
	&[\xi_2,\xi_1]=(1-q^{-1})^M(1-q^{-2})^M q_{21}^{M(M-1)} \xi_{112}; \\
	&[\xi_{112},\xi_{2}]=(1-q^{-2})^M q_{21}^{M^2} \xi_{12};\\
	&[\xi_1,\xi_{112}] =[\xi_2,\xi_{12}] =0.
	\end{align*}
Thus $\mathfrak{F}: \cU(C_2^+)\ra \az_\bq$,
$ e_1\mapsto \xi_1$, $e_2\mapsto \xi_2$, is an isomorphism of algebras.
(Of course $C_2 \simeq B_2$ but in higher rank we will get different root systems depending on the parity of $N$).

\end{row}

\begin{row}
Let $q\in\G_N'$, $N \neq 2,4$. These diagrams correspond to braidings of super type $B$ with
$\Delta^+_\bq=\{ \alpha_1, 2\alpha_1+\alpha_2, \alpha_1+\alpha_2, \alpha_2\}$.
	
If the diagram is \setlength{\unitlength}{1mm} \Dchaintwo{}{$q$}{$q^{-2}$}{$-1$},
then the Cartan roots are $\alpha_1$ and $\alpha_1+\alpha_2$, with $N_{\alpha_1}=N$, $N_{\alpha_1+\alpha_2}=M$;
here, $M=N$ if $N$ is odd and $M=\frac{N}{2}$ if $N$ is even.
The elements $x_1^N,x_{12}^M\in\dpn_\bq$ are primitive in $\dpndual$.
Thus, in $\az_\bq$, $[\xi_{12},\xi_1]=0$ and  $\az_\bq\simeq \cU((A_1\oplus A_1)^+)$.
	
If we consider the diagram \setlength{\unitlength}{1mm}	 \Dchaintwo{}{$-q^{-1}$}{$q^{2}$}{$-1$},
then $\fO_\bq=\{ \alpha_1, \alpha_1+\alpha_2\}$, $N_{\alpha_1}=M$ and  $N_{\alpha_1+\alpha_2}=N$.
The elements $x_1^M,x_{12}^N\in\dpn_\bq$ are primitive, so $[\xi_{12},\xi_1]=0$ and
$\az_\bq\simeq \cU((A_1\oplus A_1)^+)$.
\end{row}

\begin{row}
Let $q\in\G_N'$, $N \neq 3$, $\zeta\in\G_3'$. The diagram \setlength{\unitlength}{1mm}\Dchaintwo{}{$\zeta$}{$q^{-1}$}{$q$}
corresponds to a braiding of standard type $B_2$, so $\Delta^+_\bq=\{ \alpha_1, 2\alpha_1+\alpha_2, \alpha_1+\alpha_2, \alpha_2\}$.
The other diagram \setlength{\unitlength}{1.2mm}\Dchaintwo{}{$\zeta$}{$q\zeta^{-1}$}{$\zeta q^{-1}$} \, is obtained by changing the parameter
$q\leftrightarrow \zeta q^{-1}$.

The Cartan roots are $2\alpha_1+\alpha_2$ and $\alpha_2$, with $N_{2\alpha_1+\alpha_2}=M:=\ord (\zeta q^{-1})$ and $N_{\alpha_2}=N$.
The elements $x_{112}^M$, $x_{2}^N\in\dpn_\bq$ are primitive. Thus, in $\az_\bq$, we have
$[\xi_{112},\xi_2]=0$. Hence,  $\az_\bq\simeq \cU((A_1\oplus A_1)^+)$.

\end{row}

\begin{row}
Let $\zeta\in\G_3'$. The diagrams \setlength{\unitlength}{1mm} \Dchaintwo{}{$\zeta$}{$-\zeta$}{$-1$}
\, and  \Dchaintwo{}{$\zeta^{-1}$}{$-\zeta^{-1}$}{$-1$} \, correspond to braidings of standard type $B$, thus
$\Delta^+_\bq=\{ \alpha_1, 2\alpha_1+\alpha_2, \alpha_1+\alpha_2, \alpha_2\}$. In both cases $\fO_\bq$ is empty
so the corresponding Lie algebras are trivial.
\end{row}

\begin{row}
Let $\zeta\in\G_{12}'$. The diagrams of this row correspond to braidings of type $\ufo(7)$. In all cases
$\fO_\bq=\emptyset$ and the associated Lie algebras are trivial.
\end{row}

\begin{row}
Let $\zeta\in\G_{12}'$. The diagrams of this row correspond to braidings of type $\ufo(8)$.
For	\setlength{\unitlength}{1mm} \Dchaintwo{}{$-\zeta^2$}{$\zeta$}{$-\zeta^2$} \, ,
$\Delta^+_\bq=\{ \alpha_1, 2\alpha_1+\alpha_2, \alpha_1+\alpha_2, \alpha_1+2\alpha_2, \alpha_2\}$.
In this case $\fO_\bq=\{\alpha_1+\alpha_2\}$, $N_{\alpha_1+\alpha_2}=12$. Hence $\az_\bq\simeq \cU(A_1^+)$.
The same result holds for the other braidings in this row.	
\end{row}

\begin{row}
Let $\zeta\in\G_{9}'$. The diagrams of this row correspond to braidings of type $\brj(2;3)$.
If $\bq$ has diagram \setlength{\unitlength}{1mm} \Dchaintwo{}{$-\zeta$}{$\zeta^7$}{$\zeta^3$}	\, , then
$$ \Delta^+_\bq=\{ \alpha_1, 2\alpha_1+\alpha_2, 3\alpha_1+2\alpha_2, \alpha_1+\alpha_2, \alpha_1+2\alpha_2, \alpha_2\}. $$
In this case $\fO_\bq=\{ \alpha_1, \alpha_1+\alpha_2\}$ and $N_{\alpha_1}=N_{\alpha_1+\alpha_2}=18$.
Thus $[\xi_{12},\xi_1]=0$, so $\az_\bq\simeq \cU((A_1\oplus A_1)^+)$.
	
If $\bq$ has diagram \setlength{\unitlength}{1mm} \Dchaintwo{}{$\zeta^3$}{$\zeta^8$}{$-1$} \, ,
\setlength{\unitlength}{1mm} \Dchaintwo{}{$-\zeta^2$}{$\zeta$}{$-1$} \, the set of positive roots are, respectively,
\begin{align*}
&\{ \alpha_1, 2\alpha_1+\alpha_2,  3\alpha_1+2\alpha_2, 4\alpha_1+3\alpha_2, \alpha_1+\alpha_2, \alpha_2\},\\
&\{ \alpha_1, 4\alpha_1+\alpha_2, 3\alpha_1+\alpha_2, 2\alpha_1+\alpha_2, \alpha_1+\alpha_2, \alpha_2\};
\end{align*}
the Cartan roots are, respectively, $\alpha_1+\alpha_2, 2\alpha_1+\alpha_2$ and $\alpha_1, 2\alpha_1+\alpha_2$.
Hence, in both cases, $\az_\bq\simeq \cU((A_1\oplus A_1)^+)$.
\end{row}

\begin{row}
Let $q\in\G_N'$, $N \geq 4$. The diagram \setlength{\unitlength}{1mm} \Dchaintwo{}{$q$}{$q^{-3}$}{$q^3$} \,
corresponds to a braiding of Cartan type $G_2$, so $\fO_\bq=\Delta^+_\bq=\{ \alpha_1, \alpha_1+\alpha_2, 2\alpha_1+\alpha_2, 3\alpha_1+\alpha_2, 3\alpha_1+2\alpha_2, \alpha_2\}$. The coproducts of the PBW-generators are:
\begin{align*}
	\underline{\Delta}&(x_1)=x_1\ot 1+ 1\ot x_1; \qquad
	\underline{\Delta}(x_2)=x_2\ot 1+ 1\ot x_2; \\
	\underline{\Delta}&(x_{12})=x_{12}\ot 1+ 1\ot x_{12}
	+(1-q^{-3})\, x_1\ot x_2; \\
	\underline{\Delta}&(x_{112})=x_{112}\ot 1+ 1\ot x_{112}
	+(1+q)(1-q^{-2})\, x_1\ot x_{12}\\ \qquad
	&+(1-q^{-2})(1-q^{-3})\, x_1^2\ot x_2;\\
	\underline{\Delta}&(x_{1112})=x_{1112}\ot 1+ 1\ot x_{1112}
	+q^2(1-q^{-3})\, x_1\ot x_{112} \\ \qquad &+(q^2-1)(1-q^{-3})\, x_1^2\ot x_{12}+(1-q^{-3})(1-q^{-2})(1-q^{-1})\, x_1^3\ot x_{2};\\
	\underline{\Delta}&([x_{112},x_{12}]_c)=[x_{112},x_{12}]_c\ot 1+ 1\ot [x_{112},x_{12}]_c
	+(q-q^{-1})\, x_{112}\ot x_{12} \\&+(1-q^{-3})(1+q)(1-q^{-1}+q)\, x_{112}x_1\ot x_{2}\\&
	-qq_{21}(1-q^{-3})(1+q-q^{2})\, x_{1112}\ot x_{2}
	+ q^2q_{21}(1-q^{-3})\, x_1\ot [x_{112},x_2]_c\\ &
	+(1-q^{-3})^2(q^2-1)\, x_1^2\ot x_2x_{12}
	\\ &+
	q_{21}(1-q^{-3})^2(1-q^{-2})(1-q^{-1})\, x_1^3\ot x_2^2.
	\end{align*}
	We have two cases.

\begin{enumerate}
  \item If $3$ does not divide $N$, then $N_\beta=N$ for all $\beta\in\Delta^+_\bq$. Thus, in $\dpn_\bq$,
 \end{enumerate}
	\begin{align*}
	\underline{\Delta}&(x_1^N)=x_1^N\ot 1+ 1\ot x_1^N; \qquad
	\underline{\Delta}(x_2^N)=x_2^N\ot 1+ 1\ot x_2^N; \\
	\underline{\Delta}&(x_{12}^N)=x_{12}^N\ot 1+ 1\ot x_{12}^N
	+a_1\, x_1^N\ot x_2^N; \\
	\underline{\Delta}&(x_{112}^N)=x_{112}^N\ot 1+ 1\ot x_{112}^N
	+a_2\, x_1^N\ot x_{12}^N
	+a_3\, x_1^{2N}\ot x_2^N;\\
	\underline{\Delta}&(x_{1112}^N)=x_{1112}^N\ot 1+ 1\ot x_{1112}^N
	+a_4\, x_1^N\ot x_{112}^N
	+a_5\, x_1^{2N}\ot x_{12}^N\\&+
	a_6\, x_1^{3N}\ot x_{2}^N;\\
	 \underline{\Delta}&([x_{112},x_{12}]_c^N)=[x_{112},x_{12}]_c^N\ot 1+ 1\ot [x_{112},x_{12}]_c^N
	+a_7\, x_{112}^N\ot x_{12}^N
	\\&+a_8\, x_{1112}^N\ot x_{2}^N
	+a_9\, x_1^{N}\ot x_{12}^{2N}
	+a_{10}\, x_1^{2N}\ot x_{2}^Nx_{12}^N\\&+a_{11}\, x_{112}^Nx_1^{N}\ot x_2^N
	+a_{12}\, x_1^{3N}\ot x_2^{2N};
	\end{align*}
for some $a_i\in\ku$. Since
\begin{align*}
a_1&=(1-q^{-3})^Nq_{21}^{\frac{N(N-1)}{2}}\neq0, \\
a_3&=(1-q^{-2})^N(1-q^{-3})^N\neq0, \\
a_6&=(1-q^{-1})^N(1-q^{-2})^N(1-q^{-3})^Nq_{21}^{\frac{3N(N-1)}{2}}\neq0, \\
a_{12}&=(1-q^{-1})^N(1-q^{-2})^N(1-q^{-3})^{2N}\neq0,
\end{align*}
the elements $x_{12}^N$, $x_{112}^N$, $x_{1112}^N$ and $[x_{112},x_{12}]_c^N$ are not primitive.
Hence $\az_\bq$ is generated by $\xi_1$ and $\xi_2$; also
\begin{align*}
[\xi_2,\xi_1]&=a_1\, \xi_{12}; & [\xi_{12},\xi_{1}] & =a_2\, \xi_{112};
\\
[\xi_{112},\xi_{1}]&=a_4\, \xi_{1112}; & [\xi_1,\xi_{1112}]&=[\xi_2,\xi_{12}] =0.
\end{align*}
Thus, we have  $\az_\bq\simeq \cU(G_2^+)$.

\begin{enumerate}
	\item[(2)] If $N=3M$, then $N_{\alpha_1}=N_{\alpha_1+\alpha_2}=N_{2\alpha_1+\alpha_2}=N$ and
$N_{3\alpha_1+\alpha_2}=N_{3\alpha_1+2\alpha_2}=N_{\alpha_2}=M$. In this case we have
\end{enumerate}	
\begin{align*}
	\underline{\Delta}&(x_1^N)=x_1^N\ot 1+ 1\ot x_1^N; \qquad
	\underline{\Delta}(x_2^M)=x_2^M\ot 1+ 1\ot x_2^M; \\
	\underline{\Delta}&(x_{12}^N)=x_{12}^N\ot 1+ 1\ot x_{12}^N
	+(1-q^{-3})^Mq_{21}^{\frac{N(M-1)}{2}}  [x_{112},x_{12}]_c^M\ot x_2^{M}\\&
	+ (1-q^{-3})^{2M} x_{1112}^M\ot x_2^{2M}
	+(1-q^{-3})^Nq_{21}^{\frac{N(N-1)}{2}} x_1^N\ot x_2^{3M}; \\
	\underline{\Delta}&(x_{112}^N)=x_{112}^M\ot 1+ 1\ot x_{112}^M
	+b_1\,x_1^{N}\ot x_{12}^N
	+ b_2\,
	x_{1112}^M\ot [x_{112},x_{12}]_c^M\\&\qquad
	+b_3\, x_1^{2N}\ot x_2^{3M} + b_4\, x_{1112}^{2M}\ot x_2^M \\
	&\qquad+ b_5\, x_{1112}^Mx_1^N \ot x_2^{2M}+ b_6\, x_{1}^N \ot x_2^{M};\\
	\underline{\Delta}&(x_{1112}^M)=x_{112}^M\ot 1+ 1\ot x_{112}^M
	+b_7\, x_1^{N}\ot x_2^M[x_{112},x_{12}]_c^M;\\
	\underline{\Delta}&([x_{112},x_{12}]_c^M)=x_{112}^M\ot 1+ 1\ot x_{112}^M
	+b_8\, x_1^{N}\ot x_2^{2M} +   b_9\,x_{1112}^{M}\ot x_2^{M};
	\end{align*}
for some $b_i\in\ku$. We compute some of them explicitly:
\begin{align*}
b_2&=(1+q)^M(1-q^{-2})^Mq^{2M}q_{21}^{\frac{N(M-1)}{2}}, \\ b_7&=(1-q^{-3})^M(1-q^{-2})^M(1-q^{-1})^Mq_{21}^{\frac{N(M-1)}{2}}, \\
b_8&= (1-q^{-3})^{2M}(1-q^{-2})^M(1-q^{-1})^Mq_{21}^M.
\end{align*}
As these scalars are not zero, the elements $x_{12}^N$, $x_{112}^N$, $x_{1112}^M$ and $[x_{112},x_{12}]_c^M$ are not primitive.
Thus $\az_\bq\simeq \cU(G_2^+)$.
\end{row}

\begin{row}
Let $\zeta\in\G'_8$. The diagrams of this row correspond to braidings of standard type $G_2$, so
$\Delta^+_\bq=\{ \alpha_1, 3\alpha_1+\alpha_2, 2\alpha_1+\alpha_2, 3\alpha_1+2\alpha_2, \alpha_1+\alpha_2, \alpha_2\}$.
	
If $\bq$ has diagram \setlength{\unitlength}{1mm} \Dchaintwo{}{$\zeta^2$}{$\zeta$}{$\zeta^{-1}$}	\, ,
then the Cartan roots are $2\alpha_1+\alpha_2$ and $\alpha_2$ with $N_{2\alpha_1+\alpha_2}=N_{\alpha_2}=8$.
The elements $x_{112}^8,x_{2}^8\in\dpn_\bq$ are primitive and $[\xi_{2},\xi_{112}]=0$ in $\az_\bq$.
Hence $\az_\bq\simeq \cU((A_1\oplus A_1)^+)$. An analogous result holds for the other diagrams of the row.
\end{row}

\begin{row}
Let $\zeta\in\G'_{24}$. This row corresponds to type $\ufo(9)$. If $\bq$ has diagram
\setlength{\unitlength}{1mm} \Dchaintwo{}{$\zeta^6$}{$\zeta^{11}$}{$\zeta^{8}$}	\, , then
$$ \Delta^+_\bq=\{ \alpha_1, 3\alpha_1+\alpha_2, 2\alpha_1+\alpha_2, 3\alpha_1+2\alpha_2, 4\alpha_1+3\alpha_2, \alpha_1+\alpha_2, \alpha_1+2\alpha_2, \alpha_2\}$$
and $\fO_\bq=\{\alpha_1+\alpha_2, 3\alpha_1+\alpha_2\}$.
Here, $N_{\alpha_1+\alpha_2}=N_{3\alpha_1+\alpha_2}=24$, and $x_{12}^{24},x_{1112}^{24}\in\dpn_\bq$ are primitive.
In $\az_\bq$ we have the relation $[\xi_{12},\xi_{1112}]=0$; thus $\az_\bq\simeq \cU((A_1\oplus A_1)^+)$.
	
For the other diagrams, \setlength{\unitlength}{1mm}\Dchaintwo{}{$\zeta^6$}{$\zeta$}{$\zeta^{-1}$} \, ,
\setlength{\unitlength}{1mm}\Dchaintwo{}{$\zeta^8$}{$\zeta^5$}{$-1$}\, and \setlength{\unitlength}{1mm}
\Dchaintwo{}{$\zeta$}{$\zeta^{19}$}{$-1$} \, , the sets of positive roots are, respectively,
\begin{align*}
&\{ \alpha_1, \alpha_1+\alpha_2, 2\alpha_1+\alpha_2, 3\alpha_1+\alpha_2, 3\alpha_1+2\alpha_2, 5\alpha_1+2\alpha_2, 5\alpha_1+3\alpha_2, \alpha_2\},\\
&\{ \alpha_1, \alpha_1+\alpha_2, 2\alpha_1+\alpha_2, 3\alpha_1+2\alpha_2, 4\alpha_1+3\alpha_2, 5\alpha_1+3\alpha_2, 5\alpha_1+4\alpha_2, \alpha_2\},\\
&\{ \alpha_1, \alpha_1+\alpha_2, 2\alpha_1+\alpha_2, 3\alpha_1+\alpha_2, 4\alpha_1+\alpha_2, 5\alpha_1+\alpha_2, 5\alpha_1+2\alpha_2, \alpha_2\}.
\end{align*}
The Cartan roots are, respectively, $2\alpha_1+\alpha_2, \alpha_2$; $\alpha_1+\alpha_2, 5\alpha_1+3\alpha_2$; $\alpha_1, 5\alpha_1+2\alpha_2$. Hence,
in all cases, $\az_\bq\simeq\cU((A_1\oplus A_1)^+)$.
\end{row}

\begin{row}
Let $\zeta\in\G'_{5}$. The braidings in this row are associated to the Lie superalgebra $\brj(2;5)$ \cite[\S 5.2]{A-preNichols}.
If $\bq$ has diagram \setlength{\unitlength}{1mm} \Dchaintwo{}{$\zeta$}{$\zeta^2$}{$-1$} \, , then
$\Delta^+_\bq=\{ \alpha_1, 3\alpha_1+\alpha_2, 2\alpha_1+\alpha_2, 5\alpha_1+3\alpha_2, 3\alpha_1+2\alpha_2, 4\alpha_1+3\alpha_2, \alpha_1+\alpha_2, \alpha_2\}$.
In this case the Cartan roots are $\alpha_1$, $\alpha_1+\alpha_2$, $2\alpha_1+\alpha_2$ and $3\alpha_1+\alpha_2$, with
$N_{\alpha_1}=N_{3\alpha_1+2\alpha_2}=5$ and $N_{\alpha_1+\alpha_2}=N_{2\alpha_1+\alpha_2}=10$. In $\dpndual$,
	\begin{align*}
	\underline{\Delta}&(x_1)=x_1\ot 1+ 1\ot x_1; \\
	\underline{\Delta}&(x_{12})=x_{12}\ot 1+ 1\ot x_{12}
	+(1-\zeta^2)\, x_1\ot x_2; \\
	\underline{\Delta}&(x_{112})=x_{112}\ot 1+ 1\ot x_{112}
	+(1+\zeta)(1-\zeta^3)\, x_1\ot x_{12} \\ \qquad
	&+(1-\zeta^{2})(1-\zeta^{3})\, x_1^2\ot x_2;\\
	\underline{\Delta}&([x_{112},x_{12}]_c)=[x_{112},x_{12}]_c\ot 1+ 1\ot [x_{112},x_{12}]_c
	\\&-\zeta^3(1-\zeta^3)(1+\zeta)^2\, x_1\ot x_{12}^2 -\zeta q_{21}\, x_1x_{112}\ot x_{2}
	\\&+(1+q_{21}+\zeta^3q_{21})\, x_{112}x_1\ot x_{2} +\zeta(1-\zeta^2)\, x_1x_{12}x_1\ot x_{2}\\&
	+  (1-\zeta^2)(1-\zeta^3)^2\, x_1^2\ot x_2x_{12}.
	\end{align*}
Hence the coproducts of $x_1^5,x_{12}^{10},x_{112}^{10},[x_{112},x_{12}]_c^5,\in\dpn_\bq$ are:
	\begin{align*}
	\underline{\Delta}&(x_1^5)=x_1^5\ot 1+ 1\ot x_1^5; \qquad
	\underline{\Delta}(x_{12}^{10})=x_{12}^{10}\ot 1+ 1\ot x_{12}^{10}; \\
	\underline{\Delta}&(x_{112}^{10})=x_{112}^{10}\ot 1+ 1\ot x_{112}^{10}
	+a_1 \, x_1^{10}\ot x_{12}^{10}
	+ a_2\, x_1^5\ot [x_{112},x_{12}]_c^5; \\
	\underline{\Delta}&([x_{112},x_{12}]_c^5)=[x_{112},x_{12}]_c^5\ot 1+ 1\ot [x_{112},x_{12}]_c^5
	+a_3\, x_1^{5}\ot x_{12}^{10}.
	\end{align*}
for some $a_i\in\ku$. Thus, the following relations hold in $\az_\bq$
\begin{align*}
[\xi_{12},\xi_{1}]&=a_3\, \xi_{112,12}; & [\xi_{112,12},\xi_1]&=a_2\, \xi_{112}; &
[\xi_1,\xi_{112,12}]&=[\xi_{12},\xi_{112}] =0.
\end{align*}
Since
\begin{align*}
a_1=& -(1-\zeta^3)^5(1+\zeta)^5(1+62\zeta-15\zeta^2-87\zeta^3+70\zeta^4)\neq0; \\ a_3=&-(1-\zeta^3)^5(1+\zeta)^8(4-8\zeta-19\zeta^2-3\zeta^3-50\zeta^4)\neq0,
\end{align*}
the elements $x_{112}^{10}, [x_{112},x_{12}]_c^5$ are not primitive, so $\xi_{1},\xi_{12}$ generate $\az_\bq$.
Hence, $\az_\bq\simeq \cU(B_2^+)$.
	
If $\bq$ has diagram \setlength{\unitlength}{1mm} \Dchaintwo{}{$-\zeta^3$}{$\zeta^3$}{$-1$}\, , then
\begin{align*}
\Delta^+_\bq&=\{ \alpha_1, 4\alpha_1+\alpha_2, 3\alpha_1+\alpha_2, 5\alpha_1+2\alpha_2, 2\alpha_1+\alpha_2, 3\alpha_1+2\alpha_2, \alpha_1+\alpha_2, \alpha_2\}, \\
\fO_\bq &=\{\alpha_1, 3\alpha_1+\alpha_2, 2\alpha_1+\alpha_2, \alpha_1+\alpha_2\},
\end{align*}
with $N_{\alpha_1}=N_{\alpha_1+\alpha_2}=10$, $N_{3\alpha_1+\alpha_2}=N_{\alpha_1+\alpha_2}=5$. The generators of $\az_\bq$ are
$\xi_1$ and $\xi_{12}$ and they satisfy the following relations
\begin{align*}
[\xi_{12},\xi_{1}]&=b_1\, \xi_{1112}, & [\xi_{1112},\xi_{12}]&=b_2\, \xi_{112}, &
[\xi_1,\xi_{1112}]&=[\xi_{12},\xi_{112}] =0,
\end{align*}
for some $b_1,b_2\in\ku^{\times}$.  Hence  $\az_\bq\simeq \cU(C_2^+)$.
\end{row}

\begin{row}
Let $\zeta\in\G'_{20}$. This row corresponds to type $\ufo(10)$. If $\bq$ has diagram \setlength{\unitlength}{1mm}
\Dchaintwo{}{$\zeta$}{$\zeta^{17}$}{$-1$} \, , then
$\Delta^+_\bq=\{ \alpha_1, 3\alpha_1+\alpha_2, 2\alpha_1+\alpha_2, 5\alpha_1+3\alpha_2, 3\alpha_1+2\alpha_2, 4\alpha_1+3\alpha_2, \alpha_1+\alpha_2, \alpha_2\}$. The Cartan roots are $\alpha_1$ and $3\alpha_1+2\alpha_2$ with $N_{\alpha_1}=N_{3\alpha_1+2\alpha_2}=20$.
The elements $x_{1}^{20},[x_{112},x_{12}]_c^{20}\in\dpn_\bq$ are primitive; thus $[\xi_{12},\xi_{112,12}]=0$ in $\az_\bq$ and
$\az_\bq\simeq \cU((A_1\oplus A_1)^+)$. The same holds when the diagram of $\bq$ is another one in this row:
$\az_\bq\simeq \cU((A_1\oplus A_1)^+)$.
\end{row}

\begin{row}
Let $\zeta\in\G'_{15}$. This row corresponds to type $\ufo(11)$. If $\bq$ has diagram
\setlength{\unitlength}{1mm} \Dchaintwo{}{$-\zeta$}{$-\zeta^{12}$}{$\zeta^{5}$}	\, , then
$\Delta^+_\bq=\{ \alpha_1, 3\alpha_1+\alpha_2, 5\alpha_1+2\alpha_2, 2\alpha_1+\alpha_2, 3\alpha_1+2\alpha_2, \alpha_1+\alpha_2, \alpha_1+2\alpha_2, \alpha_2\}$. The Cartan roots are $\alpha_1$ and $3\alpha_1+2\alpha_2$ with $N_{\alpha_1}=N_{3\alpha_1+2\alpha_2}=30$.
In $\az_\bq$ we have $[\xi_{12},\xi_{112,12}]=0$, thus $\az_\bq\simeq \cU((A_1\oplus A_1)^+)$.
The same result holds if we consider the other diagrams of this row.
\end{row}

\begin{row}
Let $\zeta\in\G'_7$. This row corresponds to type $\ufo(12)$. If $\bq$ has diagram
\setlength{\unitlength}{1mm} \Dchaintwo{}{$-\zeta^5$}{$-\zeta^3$}{$-1$} \, , then
\begin{multline*}
\Delta^+_\bq=\{ \alpha_1, 5\alpha_1+\alpha_2, 4\alpha_1+\alpha_2, 7\alpha_1+2\alpha_2, 3\alpha_1+\alpha_2, 8\alpha_1+3\alpha_2,
\\
5\alpha_1+2\alpha_2, 7\alpha_1+3\alpha_2, 2\alpha_1+\alpha_2, 3\alpha_1+2\alpha_2, \alpha_1+\alpha_2, \alpha_2\}.
\end{multline*}
Also, $\fO_\bq=\{\alpha_1, 4\alpha_1+\alpha_2, 3\alpha_1+\alpha_2, 5\alpha_1+2\alpha_2, 2\alpha_1+\alpha_2, \alpha_1+\alpha_2\}$
with $N_{\beta}=14$ for all $\beta\in\fO_\bq$. In $\dpndual$ we have
\begin{align*}
\underline{\Delta}&(x_1)=x_1\ot 1+ 1\ot x_1; \\
\underline{\Delta}&(x_{12})=x_{12}\ot 1+ 1\ot x_{12} +(1+\zeta^3)\, x_1\ot x_2; \\
\underline{\Delta}&(x_{112})=x_{112}\ot 1+ 1\ot x_{112} +(1-\zeta)(1-\zeta^5)\, x_1\ot x_{12} \\
\qquad &+(1-\zeta)(1+\zeta^{3})\, x_1^2\ot x_2;\\
\underline{\Delta}&(x_{1112})=x_{1112}\ot 1+ 1\ot x_{1112} +(1+\zeta^3-\zeta^5)(1+\zeta^6)\, x_1\ot x_{112}	\\
& +\zeta(\zeta^3-1)\, x_1^2\ot x_{12} +\zeta^6(1-\zeta^2)(1+\zeta^3)\, x_{1}^3\ot x_{2};\\
\underline{\Delta}&(x_{11112})=x_{11112}\ot 1+ 1\ot x_{11112} -\zeta(1-\zeta)(1-\zeta^2)\, x_1\ot x_{1112}\\
& +(1-\zeta^4)\, x_1^2\ot x_{112}-(1-\zeta)(1-\zeta^2)^2\, x_{1}^3\ot x_{12}\\
& +\zeta^2(1-\zeta)(1-\zeta^2)\, x_{1}^4\ot x_{2};\\
\underline{\Delta}&([x_{1112},x_{112}]_c)=[x_{1112},x_{112}]_c\ot 1+ 1\ot [x_{1112},x_{112}]_c\\
& -\frac{(1-\zeta^5)}{(1+\zeta)}(1-\zeta^3+2\zeta^4)\, x_1\ot x_{112}^2 \\
& -q_{21}(1-\zeta)(1-\zeta^3)\, x_1^2\ot [x_{112},x_{12}]_c\\
& -(1-\zeta)^2(4+4\zeta+\zeta^2-2\zeta^3-3\zeta^4)\, x_1^2\ot x_{12}x_{112}\\
& +q_{21}(1-\zeta^2)^2\zeta^4(1-2\zeta-3\zeta^4-2\zeta^5+\zeta^6)\, x_1^3\ot x_{12}^2\\
& +(1-\zeta)^2(1+\zeta^3)^2(1+\zeta^6)\, x_1^3\ot x_2x_{112}-\zeta(1-\zeta)(1-\zeta^2)\, x_{1112}\ot x_{112}\\
& -q_{21}\zeta^6(1-\zeta)^2(1-\zeta^2)(1+2\zeta)\, x_1^4\ot x_2x_{12}\\
& +q_{21}^2\zeta^2(1-\zeta)^2(1-\zeta^2)(1+\zeta^3)\, x_1^5\ot x_2^{2}\\
& -q_{12}^2(1+\zeta^3)(1-\zeta)(1-\zeta^4+\zeta^6)\, x_{111112}\ot x_2 \\
& + \zeta q_{21}(1+\zeta^3)(1-\zeta)(1-\zeta^2)(1+\zeta-\zeta^2)\, x_{11112}x_1\ot x_2\\
& -\zeta(1-\zeta)^2(1+\zeta^3)(1-\zeta-2\zeta^2-\zeta^3)\,x_{1112}x_1^2\ot x_2\\
& +(1-\zeta)(1+\zeta^2+\zeta^3-\zeta^4-\zeta^5)\, x_{1112}x_1\ot x_{12}\\
& +\zeta q_{21}(1-\zeta)^2(2+\zeta-\zeta^3)\, x_{11112}\ot x_{12}.
\end{align*}

\smallbreak
Hence
	\begin{align*}
	\underline{\Delta}&(x_1^{14})=x_1^{14}\ot 1+ 1\ot x_1^{14}; \qquad
	\underline{\Delta}(x_{12}^{14})=x_{12}^{14}\ot 1+ 1\ot x_{12}^{14};\\
	\underline{\Delta}&(x_{112}^{14})=x_{112}^{14}\ot 1+ 1\ot x_{112}^{14}+ a_1\, x_1^{14}\ot x_{12}^{14};
	\\
	\underline{\Delta}&(x_{1112}^{14})=x_{1112}^{14}\ot 1+ 1\ot x_{1112}^{14}+ a_2\, x_1^{14}\ot x_{112}^{14}+a_3\, x_1^{28}\ot x_{12}^{14};\\
	\underline{\Delta}&(x_{11112}^{14})=x_{11112}^{14}\ot 1+ 1\ot x_{11112}^{14}+ a_4\, x_1^{14}\ot x_{1112}^{14}\\&
	\qquad +a_5\, x_1^{28}\ot x_{112}^{14}+a_6\, x_1^{42}\ot x_{12}^{14};\\
	 \underline{\Delta}&([x_{1112},x_{112}]_c^{14})=[x_{1112},x_{112}]_c^{14}\ot 1+ 1\ot[x_{1112},x_{112}]_c^{14}+a_7\, x_{1112}^{14}\ot x_{12}^{14}\\&
	\qquad +a_8\, x_{11112}^{14}\ot x_{12}^{14}+a_9\, x_1^{42}\ot x_{12}^{28}
	+a_{10}\, x_{1}^{14}\ot x_{112}^{28}\\&
	\qquad +a_{11}\, x_1^{28}\ot x_{12}^{14}x_{112}^{14}+a_{12}\, x_{1112}^{14}x_1^{14}\ot x_{12}^{14};
	\end{align*}
with $a_i\in\ku$. For instance,
\begin{gather*}
a_1=q_{21}^7(1-\zeta)^7(1-\zeta^5)^7\big(4059-7124\zeta +35105\zeta^2 +31472\zeta^3 -17431\zeta^4 
\\+19299\zeta^5 +40124\zeta^6 \big)\neq 0,
\end{gather*}
because $\zeta \in \G'_7$ . Also,

\begin{gather*}
a_3= 26686268 + 39070423 \zeta -42643895 \zeta^2 - 19103336 \zeta^3 +52678504\zeta^4 
\\- 4378676 \zeta^5 - 51111858 \zeta^6 \neq0.
\end{gather*}
Since $a_1, a_3, a_6, a_{12}\neq 0$
then $x_{112}^{14}$, $x_{1112}^{14}$, $x_{11112}^{14}$ and $[x_{1112},x_{112}]_c^{14}$ are not primitive elements in $\dpndual$. Thus, $\xi_1$ and $\xi_{12}$ generates $\az_\bq$.

Also, in $\az_\bq$ we have
	\begin{align*}
	&[\xi_{12},\xi_{1}]=a_1 \,\xi_{112};
	&[\xi_{1},\xi_{112}]=a_2 \,\xi_{1112};\\
	&[\xi_{1},\xi_{1112}]=a_4 \,\xi_{11112};
	&[\xi_1,\xi_{11112}] =[\xi_{12},\xi_{112}] =0.
	\end{align*}
	So,  $\az_\bq\simeq \cU(G_2^+)$.
	
In the case of the diagram 	
	\setlength{\unitlength}{1mm}
	\Dchaintwo{}{$-\zeta$}{$-\zeta^4$}{$-1$}
	\,
 $\az_\bq$ is generated by $\xi_{1}$, $\xi_{12}$ and
	\begin{align*}
	[\xi_{12},\xi_{1}] &=b_1 \,\xi_{112};&
	[\xi_{12},\xi_{112}]&=b_2 \,\xi_{112,12};\\
[\xi_{12},\xi_{112,12}]&=b_3 \,\xi_{(112,12),12};
	&[\xi_1,\xi_{112}] &=[\xi_{12},\xi_{(112,12),12}] =0,
	\end{align*}
where $b_1,b_2,b_3\in\ku^\times$. Hence, we also have $\az_\bq\simeq \cU(G_2^+)$.
\end{row}

\begin{remark}
The results of this paper are part of the thesis of one of the authors \cite{rossi}, where missing details of the computations can be found.
\end{remark}

\end{document}